\documentclass[11pt]{article}
\usepackage[a4paper,margin=30mm]{geometry}
\usepackage{amsmath,amssymb,amsthm}
\usepackage[hidelinks]{hyperref}
\usepackage{microtype}
\newtheorem{theorem}{Theorem}[section]
\newtheorem{proposition}[theorem]{Proposition}

\theoremstyle{remark}
\newtheorem{remark}[theorem]{Remark}
\newcommand{\Pp}{\mathbb P}

\title{Domains of attraction of quasi-stationary distributions \\
for an absorbed random walk}
\author{Phil. Pollett}
\date{}
\begin{document}
\maketitle

\begin{abstract}
We study a continuous-time random walk on the positive integers,
killed at zero, where the downward rate $a$ and upward rate $b$
satisfy $a>b>0$. Its quasi-stationary distributions form a one-parameter
family. We show that if the initial distribution $(m_n)_{n\geq 1}$ is
strictly positive and satisfies $m_{n+1}/m_n\to r$, where
$(b/a)^{1/2}<r<1$, then its conditional distribution given survival
converges in total variation to the corresponding nonminimal
quasi-stationary distribution. We obtain asymptotic equivalents for
the joint probabilities of position and survival, and for the survival
probability itself. We also give a moment condition under which the
minimal quasi-stationary distribution is selected. Finally, we
construct an initial distribution whose successive-term ratios do not
converge but whose conditional distribution nevertheless converges to
a nonminimal quasi-stationary distribution. Thus the ratio condition
is sufficient but not necessary.
\end{abstract}


\medskip


\section{Introduction}\label{sec:introduction}

We consider a continuous-time random walk on the positive integers,
killed at zero, with downward rate $a$ and upward rate $b$, where
$a>b>0$. For this process, the full family of quasi-stationary
distributions (QSDs) is known, as is the limiting conditional
distribution from each fixed initial state. Seneta established the
fixed-state limit in~\cite{Seneta1966}, while Seneta and Vere-Jones
developed a general framework for countably infinite Markov chains with
absorption~\cite{SenetaVereJones1966}. Van Doorn identified the full
family of QSDs for birth--death processes
\cite[Theorem~3.2]{VanDoorn1991} and extended the fixed-state
convergence result to initial distributions with finite support. He also
observed that treating initial distributions with infinite support would
require conditions comparing their tails with that of the minimal
QSD~\cite[p.~695]{VanDoorn1991}. Related questions for discrete-time
birth--death processes were addressed by van Doorn and
Schrijner~\cite{VanDoornSchrijner1995}, who connected limiting
conditional distributions with ratio limits for transition probabilities
through their spectral representation.

The purpose of this paper is to identify tail conditions under which an
initial distribution with infinite support selects a specified member of
the QSD family. Sufficiently light tails lead to selection of the minimal
QSD, whereas heavier, asymptotically geometric tails can lead to selection
of a nonminimal QSD. Our main result makes this precise.
Put $\rho=(b/a)^{1/2}$. If the initial distribution
$m=(m_n)_{n\geq 1}$ is strictly positive and
\[
 \frac{m_{n+1}}{m_n}\to r\in(\rho,1),
\]
then the conditional distribution of the walk given survival converges
in total variation to a corresponding nonminimal QSD $q^{(r)}$. We also
obtain asymptotic equivalents for the joint probabilities of position
and survival, and for the survival probability itself. The theorem
applies not only to geometric distributions, but also to geometric tails
modified by polynomial, logarithmic or more general subexponential
factors.

Analogous dependence on the initial distribution occurs in other
absorbing processes. For Brownian motion with negative drift,
Mart\'{\i}nez, Picco and San Mart\'{\i}n~\cite{MartinezEtAl1998} showed
that the rate at which the initial density decreases determines the
limiting QSD, and constructed an example with different limits along
different sequences of times. For subcritical branching processes, tail
conditions on the initial population distribution similarly determine
the limiting QSD \cite{RubinVereJones1968,AsmussenHering1983}. We are
not aware of a corresponding geometric-tail domain-of-attraction result
for the absorbed random walk considered here.

In earlier work, we developed a general framework for polynomial
corrections to exponential decay in absorbing continuous-time Markov
chains, using the same killed random walk as an
example~\cite{PollettAAP}. The questions addressed here are different,
but the two studies use related Bessel-function representations of the
transition probabilities of the killed process.

The ratio condition in our main theorem is sufficient but not necessary.
We construct an initial distribution whose probabilities form a
geometric sequence with parameter $r$, multiplied by a periodic factor.
Although its successive-term ratios do not converge, its conditional
distribution converges to the nonminimal QSD $q^{(r)}$. Thus, a complete
domain-of-attraction criterion cannot be formulated solely in terms of
convergence of successive-term ratios.

The paper is organised as follows. Section~\ref{sec:model} introduces
the model and its QSDs, and gives an exact calculation for a geometric
initial distribution. Section~\ref{sec:main} establishes the
geometric-tail domain-of-attraction theorem.
Section~\ref{sec:lightcritical} treats lighter tails and discusses the
critical case. Section~\ref{sec:counterexample} gives an example with
oscillating successive-term ratios, and Section~\ref{sec:necessary}
derives a necessary survival-ratio condition for convergence to a
nonminimal QSD. Section~\ref{sec:discussion} summarises the resulting
partial classification and discusses the remaining open problems.

\section{The model and its quasi-stationary distributions}
\label{sec:model}

Let $(X_t)_{t\geq 0}$ be the continuous-time random walk on the positive
integers, killed when it reaches zero. From each state $i\geq 1$, the
process jumps to $i+1$ at rate~$b$ and to $i-1$ at rate~$a$, where
$a>b>0$. In queueing terminology, this is the M/M/1 queue with arrival
rate~$b$ and service rate~$a$, stopped when the queue first empties.
Write $C=\{1,2,\ldots\}$ for the irreducible class of transient states.
Let $T=\inf\{t\geq 0:X_t=0\}$ be the absorption time. Write
$P_t=(p_{ij}(t))_{i,j\geq 1}$ for the substochastic transition semigroup
of the killed process, where
\[
 p_{ij}(t)=\mathbb P_i(X_t=j,T>t),
\]
and write $Q$ for its generator.

Seneta~\cite{Seneta1966} proved that, for each fixed initial state
$i\geq 1$,
\[
 \mathbb P_i(X_t=j\mid T>t)\to 
 q_j^*:=(1-\rho)^2j\rho^{j-1},
 \qquad j\geq 1,
\]
as $t\to\infty$, where $\rho=(b/a)^{1/2}$. Van Doorn's extension of
the fixed-state result shows that the same conclusion holds for every
initial distribution with finite support~\cite[Section~4]{VanDoorn1991}. This \emph{limiting
conditional distribution} is one of a family of quasi-stationary
distributions. A probability distribution $q=(q_j)_{j\geq 1}$ is called
a \emph{quasi-stationary distribution} (QSD) if
\[
 \mathbb P_q(X_t=j\mid T>t)=q_j,
 \qquad j\geq 1,\quad t\geq 0.
\]

Standard theory (see, for example,~\cite{VanDoornPollett2013}) shows
that this is equivalent to $q$ being an $x$-invariant probability
distribution, that is,
$qP_t=e^{-xt}q$, $t\geq 0$, or, componentwise,
\[
 \sum_{i\geq 1}q_i p_{ij}(t)=e^{-xt}q_j,
 \qquad j\geq 1,\quad t\geq 0,
\]
for some $x>0$. Summing over $j$ gives
$\mathbb P_q(T>t)=e^{-xt}$, $t\geq 0$, which says that $x$ is the
exponential decay rate of the survival probability when $q$ is taken
as the initial distribution.

The QSDs for an absorbing birth--death process are described
in~\cite[Theorem~16]{VanDoornPollett2013}. 
For the present constant-rate process, the equation $qQ=-xq$ becomes
\[
 bq_{j-1}+aq_{j+1}-(a+b)q_j=-xq_j,
 \qquad j\geq 1,
\]
with $q_0=0$. It is easy to show that the normalized solutions are
\begin{equation}
\label{eq:qsd-family}
 q_j^{(r)}
 =\frac{(1-r)(1-s)}{r-s}(r^j-s^j),
 \qquad j\geq 1,
\end{equation}
where $r\in(\rho,1)$, $s=\rho^2/r$, and
\[
 x(r)=a+b-ar-\frac br.
\]
In particular, $0<s<\rho<r<1$.

As $r\downarrow\rho$, we have $s\uparrow\rho$, so the two parameters
coalesce at $\rho$. At the same time, $x(r)$ increases to the
\emph{decay parameter} of the class $C$, given by
\[
 \lambda_C=x(\rho)
 =a+b-2\sqrt{ab}
 =(\sqrt a-\sqrt b)^2,
\]
while \eqref{eq:qsd-family} converges componentwise to $q^*$. Thus $q^*$
is the family member with the largest decay rate, $\lambda_C$. It is
called the \emph{minimal QSD} because it is stochastically smaller than
every other QSD in the family.

All members of the family can be limiting conditional distributions for
an appropriately chosen initial distribution $m=(m_n)_{n\geq 1}$; we say
that $m$ lies in the \emph{domain of attraction} of the QSD $q$ if
\[
 \mathcal L_m(X_t\mid T>t)\to q
\]
as $t\to\infty$; sometimes we say that $m$ \emph{selects} $q$. Thus the results of Seneta and van Doorn show that every initial
distribution with finite support lies in the domain of attraction of the
minimal QSD~$q^*$. Every QSD has
a nonempty domain of attraction, since a QSD is invariant under
conditioning on survival: in particular, $q^{(r)}$ itself selects
$q^{(r)}$. The substantive question is whether one can identify natural
classes of initial distributions that select a given member of the
family. We next show that geometric initial distributions provide such a
class.

Let $m_j^{(r)}=(1-r)r^{j-1}$, $j\geq 1$, where $r\in (\rho, 1)$. 
Since 
\[
 r^j={r m_j^{(r)}}/{1-r}
 \qquad\text{and}\qquad
 s^j={s m_j^{(s)}}/{1-s},
\]
a direct calculation gives
\begin{equation}
\label{eq:convex}
 m^{(r)}=A(r)q^{(r)}+B(r)m^{(s)},
 \qquad
 A(r)=\frac{r-s}{r(1-s)},\quad
 B(r)=\frac{s(1-r)}{r(1-s)}.
\end{equation}
Both coefficients are positive and their sum is one, so this is a
convex decomposition.

The positive sequence $h_i=i\rho^{-i}$, $i\geq 1$, with $h_0=0$,
satisfies
\begin{equation}
 bh_{i+1}+ah_{i-1}-(a+b)h_i=-\lambda_C h_i,
 \qquad i\geq 1.
\label{eq:h}
\end{equation}
The associated Doob $h$-transformed process is a birth-death process
with birth and death rates
\[
 \sqrt{ab}\,\frac{i+1}{i}
 \quad\hbox{and}\quad
 \sqrt{ab}\,\frac{i-1}{i},
\]
respectively, and is honest. Hence
\cite[Theorem 3.1(ii)]{Pollett1988}
$P_th=e^{-\lambda_C t}h$.
Thus, whenever $\sum_i\mu_i h_i<\infty$,
\[
 \mathbb P_\mu(T>t)
 \leq \mathbb E_\mu[h(X_t);T>t]
 =\sum_{i\geq 1}\mu_i(P_th)_i
 =e^{-\lambda_C t}\sum_{i\geq 1}\mu_i h_i.
\]
For $\mu=m^{(s)}$, the sum on the right is finite because $s<\rho$.
Consequently,
\[
 \mathbb P_{m^{(s)}}(T>t)=O(e^{-\lambda_C t}).
\]
On the other hand, a process started with distribution $q^{(r)}$
survives to time $t$ with probability $e^{-x(r)t}$, so the
$q^{(r)}$ component in \eqref{eq:convex} contributes
$A(r)e^{-x(r)t}$ to the survival probability. Moreover,
$\lambda_C-x(r)={a(r-\rho)^2}/{r}>0$.
Thus the contribution of the $m^{(s)}$ component to the survival
probability is asymptotically negligible relative to that of the
$q^{(r)}$ component. It follows that
\[
 \mathcal L_{m^{(r)}}(X_t\mid T>t)
 \to q^{(r)}
\]
in total variation and that
\[
 \mathbb P_{m^{(r)}}(T>t)
 \sim A(r)e^{-x(r)t}.
\]

This exact calculation anticipates the general tail-ratio result of
the next section.

\section{A geometric-tail domain of attraction}\label{sec:main}

The exact calculation above suggests that an initial distribution with
an asymptotically geometric tail should select the same QSD as the
corresponding geometric initial law. 
The next theorem gives the resulting domain-of-attraction criterion.
The proof uses Seneta's explicit
formula~\cite[equation~(2.2)]{Seneta1966} for
the transition probabilities:
\[
 p_{ij}(t)
 =e^{-(a+b)t}\rho^{j-i}
 \{I_{i-j}(ct)-I_{i+j}(ct)\},
 \qquad i,j\geq 1,
\]
where $c=2\sqrt{ab}$ and
$I_n$ is the usual modified Bessel function of the first kind. It
allows the effect of the tail of the initial distribution to be
isolated.

\begin{theorem}[Ratio-tail sufficiency]\label{thm:ratio-tail}
Let $m=(m_n)_{n\geq 1}$ be a probability distribution such that $m_n>0$ and
\[
 \frac{m_{n+1}}{m_n}\to r\in(\rho,1).
\]
Then
\[
 \left\|\mathcal L_m(X_t\mid T>t)-q^{(r)}\right\|_{\rm TV}
 \to 0.
\]
More precisely, if
\[
 W(t)=e^{-(a+b)t}\sum_{n\geq 1}m_n\rho^{-n}
 I_n(2\sqrt{ab}\,t),
\]
then, for every fixed $j\geq 1$, and $t\to\infty$,
\[
 \Pp_m(X_t=j,T>t)\sim(r^j-s^j)W(t),
\]
and
\[
 \Pp_m(T>t)\sim
 \frac{r-s}{(1-r)(1-s)}W(t).
\]
In particular,
\[
 \lim_{t\to\infty} - \frac1t\log\Pp_m(T>t)=x(r).
\]
\end{theorem}

\begin{proof}
For convenience, put
$\kappa_n(t) =e^{-(a+b)t}\rho^{-n}I_n(ct)$,
$n\in\mathbb Z$.
Since $I_{-n}=I_n$, we have
$\kappa_{-n}(t)=\rho^{2n}\kappa_n(t)$, $n\geq 0$.
The transition probabilities may therefore be written as
\[
 p_{ij}(t)=\kappa_{i-j}(t)-\rho^{2j}\kappa_{i+j}(t).
\]
For $j\in\mathbb Z$, set
$w_j(t)=\sum_{i\geq 1}m_i\kappa_{i-j}(t)$.
Then, for $j\geq 1$,
\[
 u_j(t):=\Pp_m(X_t=j,T>t)
 =w_j(t)-\rho^{2j}w_{-j}(t),
\]
and
$w_0(t)=\sum_{i\geq 1}m_i\kappa_i(t)=W(t)$.
We first determine the exponential order of $W$. The ratio hypothesis
implies $\log m_n=n\log r+o(n)$. Choose $\varepsilon>0$
sufficiently small that
$\rho<r-\varepsilon<r+\varepsilon<1$.
Then the tail of $m_n$ is bounded above and below by constant multiples
of $(r+\varepsilon)^n$ and $(r-\varepsilon)^n$, respectively.
The generating function of the sequence 
$(I_n)$~\cite[9.6.33]{AbramowitzStegun} is
\[
 \sum_{n\in\mathbb Z}I_n(ct)y^n
 =\exp\left(\frac{ct}{2}\left(y+y^{-1}\right)\right),
 \qquad y>0.
\]
Applying this identity with $y=(r\pm\varepsilon)/\rho$ shows that the
bilateral sums have exponential rates $-x(r\pm\varepsilon)$. Extending
the original sums from $n\geq N$ to $n\in\mathbb Z$ adds the finitely
many terms with $1\leq n<N$ and the terms with $n\leq 0$. We now show
that these additional terms have exponential rate at most $-\lambda_C$.
For each fixed index, the standard fixed-order asymptotic formula
$I_n(ct)\sim e^{ct}/\sqrt{2\pi
ct}$~\cite[equation~(9.7.1)]{AbramowitzStegun} shows that the
associated term has exponential rate $-(a+b)+c=-\lambda_C$. Moreover, if
$y=(r\pm\varepsilon)/\rho>1$, then
\[
 e^{-(a+b)t}\sum_{n\leq 0}I_n(ct)y^n
 =e^{-(a+b)t}\sum_{k\geq 0}I_k(ct)y^{-k}
 \leq e^{-(a+b)t}\sum_{k\geq 0}I_k(ct)
 \leq e^{-\lambda_C t}.
\]
Finally,
\[
 \lambda_C-x(r)
 =ar+\frac br-2\sqrt{ab}
 =\left(\sqrt{ar}-\sqrt{\frac br}\right)^2>0.
\]
By continuity, $\lambda_C>x(r\pm\varepsilon)$ when $\varepsilon$ is
sufficiently small. Thus the additional terms are exponentially
negligible and do not affect the rates $-x(r\pm\varepsilon)$. The
geometric lower and upper bounds therefore give
\[
 -x(r-\varepsilon)
 \leq
 \liminf_{t\to\infty}\frac{1}{t}\log W(t)
 \leq
 \limsup_{t\to\infty}\frac{1}{t}\log W(t)
 \leq
 -x(r+\varepsilon).
\]
Letting $\varepsilon\downarrow0$ and using the continuity of $x$
gives
\[
 \frac{1}{t}\log W(t)\to-x(r),
\]
or equivalently, $W(t)=e^{-x(r)t+o(t)}$.
Define a probability measure $\mu_t$ on $\mathbb N$ by
\[
 \mu_t(n)=\frac{m_n\kappa_n(t)}{W(t)}.
\]
For each fixed $N$, the contribution of $1\leq n\leq N$ is
$O(t^{-1/2}e^{-\lambda_C t})$, whereas
$W(t)=e^{-x(r)t+o(t)}$. Hence
$\mu_t(\{1,\ldots,N\})\to 0$ for every fixed $N$, so the mass of
$\mu_t$ escapes to infinity.

For each fixed $j$, the sequences $m_{n+j}/m_n$ and
$m_{n-j}/m_n$, with the latter restricted to $n\geq j+1$,
converge to $r^j$ and $r^{-j}$, respectively, and are therefore
bounded. Thus, since the mass of $\mu_t$ escapes to infinity,
their expectations under $\mu_t$ converge to these limits.
For $j\geq 1$, splitting the sum defining $w_j(t)$ gives
$w_j(t)=D_j(t)+B_j(t)$,
where
\[
 D_j(t)=\sum_{n\geq 1}m_{n+j}\kappa_n(t),
 \qquad
 B_j(t)=\sum_{k=0}^{j-1}m_{j-k}\kappa_{-k}(t),
\]
while
\[
 w_{-j}(t)=\sum_{n\geq j+1}m_{n-j}\kappa_n(t).
\]
For each fixed $j\geq 1$,
\[
 \frac{D_j(t)}{W(t)}
 =\sum_{n\geq 1}\frac{m_{n+j}}{m_n}\mu_t(n)
 \to r^j,
\]
because $m_{n+j}/m_n\to r^j$ as $n\to\infty$ and the mass of
$\mu_t$ escapes to infinity. Similarly,
\[
 \frac{\rho^{2j}w_{-j}(t)}{W(t)}
 =\rho^{2j}\sum_{n\geq j+1}
 \frac{m_{n-j}}{m_n}\mu_t(n)
 \to \rho^{2j}r^{-j}=s^j.
\]
The finite remainder terms $B_j(t)$ are negligible even after
summation over $j$.
Indeed, by Tonelli's theorem, and using
$\sum_{j\geq k+1}m_{j-k}=1$ for each $k\geq 0$,
\[
 \sum_{j\geq 1}B_j(t)
 =\sum_{k\geq 0}\kappa_{-k}(t).
\]
Since $0<\rho<1$ and $I_k(ct)\geq 0$,
\[
 \sum_{k\geq 0}\kappa_{-k}(t)
 =e^{-(a+b)t}\sum_{k\geq 0}\rho^k I_k(ct)
 \leq e^{-(a+b)t}\sum_{k\geq 0}I_k(ct).
\]
The Bessel generating-function identity at $y=1$ gives
$\sum_{n\in\mathbb Z}I_n(ct)=e^{ct}$.
Consequently,
\[
 \sum_{j\geq 1}B_j(t)
 =\sum_{k\geq 0}\kappa_{-k}(t)
 \leq e^{-(a+b)t}e^{ct}
 =e^{-\lambda_C t}.
\]
Since $W(t)=e^{-x(r)t+o(t)}$ and $\lambda_C>x(r)$, it follows that
${B_j(t)}/{W(t)}\to 0$
for every fixed $j$. Consequently,
\[
 \frac{u_j(t)}{W(t)}\to r^j-s^j,
 \qquad j\geq 1.
\]
It remains to obtain summable bounds for $D_j(t)$ and
$\rho^{2j}w_{-j}(t)$. 
Choose $\rho<\eta<r<R<1$. There is an integer $N$ such that
\[
 \eta\leq\frac{m_{n+1}}{m_n}\leq R,\qquad n\geq N.
\]
Since $m_n>0$ for every $n\geq 1$, the finitely many ratios involving
indices below $N$ can be absorbed into constants: there are
constants $A_1,A_2<\infty$ such that
\[
 m_{n+j}\leq A_1R^jm_n,\qquad n,j\geq 1,
\]
and
\[
 m_{n-j}\leq A_2\eta^{-j}m_n,
 \qquad j\geq 1,\quad n\geq j+1.
\]
It follows that
\[
 \frac{D_j(t)}{W(t)}
 =\sum_{n\geq 1}\frac{m_{n+j}}{m_n}\mu_t(n)
 \leq A_1R^j
\]
and
\[
 \frac{\rho^{2j}w_{-j}(t)}{W(t)}
 =\rho^{2j}\sum_{n\geq j+1}
 \frac{m_{n-j}}{m_n}\mu_t(n)
 \leq A_2\left(\frac{\rho^2}{\eta}\right)^j.
\]
Both geometric upper bounds are summable over $j$.
Dominated convergence, applied separately to $D_j(t)$ and
$\rho^{2j}w_{-j}(t)$, together with the aggregate estimate for the
finite remainder terms $B_j(t)$, therefore gives
\begin{align*}
 \frac{\Pp_m(T>t)}{W(t)}
 &=\sum_{j\geq 1}\frac{D_j(t)}{W(t)}
  +\sum_{j\geq 1}\frac{B_j(t)}{W(t)}
  -\sum_{j\geq 1}\frac{\rho^{2j}w_{-j}(t)}{W(t)}\\
 &\to 
 \sum_{j\geq 1}r^j-\sum_{j\geq 1}s^j
 =\frac{r-s}{(1-r)(1-s)}.
\end{align*}
Finally, since $0\leq u_j(t)\leq D_j(t)+B_j(t)$, the geometric
bound on $D_j(t)$ gives
\[
 \sum_{j>J}\frac{u_j(t)}{W(t)}
 \leq
 A_1\sum_{j>J}R^j
 +
 \frac{1}{W(t)}\sum_{j\geq 1}B_j(t).
\]
The second term tends to zero by the aggregate estimate above.
Taking the superior limit and then letting $J\to\infty$ shows that
the conditioned laws are tight.
After division by the survival probability, the conditioned
probabilities therefore converge pointwise to~$q_j^{(r)}$. Since both
the conditioned probabilities and $(q_j^{(r)})_{j\geq 1}$ have total
mass one, pointwise convergence implies convergence in total variation.
\end{proof}

\begin{remark}
The proof also explains why the limiting QSD involves the difference
of two geometric terms. In the killed transition formula, the first
Bessel term contributes~$r^j$, whereas the second term, which enforces
absorption at zero, contributes
$\rho^{2j}r^{-j}=s^j$. Their difference is therefore $r^j-s^j$, as
in the formula for $q_j^{(r)}$.
\end{remark}

\begin{remark}
The theorem applies well beyond geometric initial distributions.
For example, suppose that
$m_n=c\,r^nL_n$, $n\geq 1$, where $c$ is the normalizing constant and
$L_{n+1}/L_n\to 1$.
Then $m_{n+1}/m_n\to r$, so, if $\rho<r<1$, the initial
distribution $m$ selects $q^{(r)}$. This includes polynomial and
logarithmic corrections, such as
\[
 L_n=n^\beta
 \quad\hbox{and}\quad
 L_n=(\log(n+1))^\gamma,
\]
as well as subexponential corrections of the form
\[
 L_n=\exp(c_0n^\delta),
 \qquad 0<\delta<1.
\]
A finite mixture of geometric distributions provides another example:
if its largest geometric ratio is $r_*\in(\rho,1)$, then the mixture
selects $q^{(r_*)}$.
\end{remark}

\section{Lighter tails and the critical boundary}
\label{sec:lightcritical}

The theorem deals with $r>\rho$, where the initial tail selects a
nonminimal QSD. 
For initial distributions with sufficiently light tails, the
conditional distribution instead converges to the minimal QSD $q^*$.

\begin{proposition}\label{prop:light}
Let $m=(m_n)_{n\geq 1}$ be a probability distribution satisfying
\[
 H_m:=\sum_{n\geq 1}m_n n\rho^{-n}<\infty.
\]
Then, for every fixed $j\geq 1$,
\[
 \mathbb P_m(X_t=j,T>t)
 \sim
 \gamma H_m j\rho^j t^{-3/2}e^{-\lambda_C t},
 \qquad
 \gamma=\frac{2}{\sqrt{2\pi}\,c^{3/2}},
\]
and
\[
 \mathbb P_m(T>t)
 \sim
 \gamma H_m\frac{\rho}{(1-\rho)^2}
 t^{-3/2}e^{-\lambda_C t}.
\]
Moreover,
\[
 \left\|
 \mathcal L_m(X_t\mid T>t)-q^*
 \right\|_{\rm TV}
 \to 0.
\]
In particular, these conclusions hold if $m_n>0$ and
${m_{n+1}}/{m_n}\to r<\rho$.
\end{proposition}

\begin{proof}
From the Bessel representation of the transition probabilities and
the fixed-order large-argument expansion,
\[
 I_{i-j}(z)-I_{i+j}(z)
 \sim
 \frac{2ij\,e^z}{\sqrt{2\pi}\,z^{3/2}}
\]
for each fixed $i,j\geq 1$. Consequently,
$t^{3/2}e^{\lambda_C t}p_{ij}(t) \to \gamma i\rho^{-i}j\rho^j$.
We also require a bound uniform in $i$ and $j$. The integral
representation for the modified Bessel 
functions~\cite[equation~(9.6.19)]{AbramowitzStegun} gives
\[
 I_{i-j}(z)-I_{i+j}(z)
 =
 \frac{2}{\pi}\int_0^\pi e^{z\cos\theta}
 \sin(i\theta)\sin(j\theta)\,d\theta.
\]
Since
$|\sin(n\theta)|\leq n\sin\theta$, for $0\leq\theta\leq\pi$, we have
\[
 0\leq I_{i-j}(z)-I_{i+j}(z)
 \leq
 \frac{2ij}{\pi}
 \int_0^\pi e^{z\cos\theta}\sin^2\theta\,d\theta
 =
 \frac{2ijI_1(z)}{z},
\]
where the last equality follows from
\cite[equations~(9.6.19) and~(9.6.26)]{AbramowitzStegun}.
Since $z^{1/2}e^{-z}I_1(z)$ is continuous on $[c,\infty)$ and
has a finite limit as $z\to\infty$, it is bounded there. Hence, for
$t\geq 1$,
\[
 t^{3/2}e^{\lambda_C t}p_{ij}(t)
 \leq A i\rho^{-i}j\rho^j
\]
for a constant $A$ independent of $i$, $j$ and $t$.
Let
\[
 u_j(t)=\mathbb P_m(X_t=j,T>t)
       =\sum_{i\geq 1}m_i p_{ij}(t).
\]
Dominated convergence in the sum over $i$ gives
$t^{3/2}e^{\lambda_C t}u_j(t) \to \gamma H_mj\rho^j$
for each $j\geq 1$, while
$t^{3/2}e^{\lambda_C t}u_j(t) \leq AH_mj\rho^j$.
Since
$\sum_{j\geq 1}j\rho^j=\rho/(1-\rho)^2$,
a second application of dominated convergence gives
\[
 \mathbb P_m(T>t)
 \sim
 \gamma H_m\frac{\rho}{(1-\rho)^2}
 t^{-3/2}e^{-\lambda_C t}.
\]
The same domination gives convergence in $\ell^1$ of the scaled
unnormalized distributions. After normalization, their limit is
\[
 \frac{j\rho^j}{\sum_{k\geq 1}k\rho^k}
 =(1-\rho)^2j\rho^{j-1}
 =q_j^*,
\]
and hence
\[
 \left\|
 \mathcal L_m(X_t\mid T>t)-q^*
 \right\|_{\rm TV}
 \to0.
\]
Finally, suppose that $m_n>0$ and
$m_{n+1}/m_n\to r<\rho$. Choose $r<r_0<\rho$.
Then, for some $A_0<\infty$, $m_n\leq A_0r_0^n$ eventually,
and therefore
$\sum_{n\geq 1}m_n n\rho^{-n}<\infty$.
The preceding conclusions now apply.
\end{proof}

\begin{remark}\label{rem:critical}
For the critical case $m_{n+1}/m_n\to\rho$ the ratio condition alone
does not determine whether $\sum_{n\geq 1}m_n n\rho^{-n}$ is finite. For
example, if $m_n=c\,\rho^n n^{-\beta}$, $n\geq 1$, where $c$ is the
normalizing constant, then
\[
 \sum_{n\geq 1}m_n n\rho^{-n}
 =c\sum_{n\geq 1}n^{1-\beta}.
\]
Thus Proposition~\ref{prop:light} applies when $\beta>2$, and these
critical-ratio initial distributions select~$q^*$. More generally,
polynomial or slowly varying factors multiplying $\rho^n$ may
determine whether the moment condition is satisfied.
\end{remark}

\section{The ratio condition is not necessary}\label{sec:counterexample}
Theorem~\ref{thm:ratio-tail} gives a simple sufficient condition, but
the raw ratios can oscillate without affecting the limiting QSD. The
next example isolates the oscillatory mode and shows that it decays too
quickly to contribute to the conditioned limit.

\begin{proposition}\label{prop:oscillatory}
Fix $r\in(\rho,1)$ and $0<\varepsilon<1$, and let
\[
 m_j=Ar^{j-1}(1+\varepsilon(-1)^j),
 \qquad j\geq 1,
\]
where
\[
 A^{-1}
 =\sum_{j\geq 1}r^{j-1}(1+\varepsilon(-1)^j)
 =\frac{1}{1-r}-\frac{\varepsilon}{1+r}.
\]
Then $m=(m_j)_{j\geq 1}$ is a strictly positive probability
distribution,
\[
 \frac{m_{j+1}}{m_j}
 =r\frac{1-\varepsilon(-1)^j}{1+\varepsilon(-1)^j}
\]
has two distinct subsequential limits, and
\[
 \left\|
 \mathcal L_m(X_t\mid T>t)-q^{(r)}
 \right\|_{\rm TV}
 \to 0.
\]
\end{proposition}

\begin{proof}
Put $s=\rho^2/r$. For $z\neq 0$, define
\[
 z^\dagger=\frac{\rho^2}{z},
 \qquad
 V_j^{(z)}=z^j-(z^\dagger)^j,
 \qquad
 x(z)=a+b-az-\frac bz.
\]
If $|z|<1$ and $|z^\dagger|<1$, then $V^{(z)}$ belongs to
$\ell^1$. A direct substitution gives
\[
 bV_{j-1}^{(z)}+aV_{j+1}^{(z)}-(a+b)V_j^{(z)}
 =-x(z)V_j^{(z)},
 \qquad j\geq 1,
\]
with $V_0^{(z)}=0$. 
Thus $V^{(z)}$ is an $\ell^1$ left eigenvector of $Q$, with
eigenvalue $-x(z)$. Since the killed generator is a bounded operator
on $\ell^1$, with
\[
 \|vQ\|_1\leq 2(a+b)\|v\|_1,
 \qquad v\in\ell^1,
\]
the semigroup is given by the operator-norm convergent series
\[
 P_t=e^{tQ}=\sum_{n=0}^{\infty}\frac{t^nQ^n}{n!}.
\]
Since this series converges in operator norm and $V^{(z)}\in\ell^1$,
left multiplication by $V^{(z)}$ may be passed through the sum.
The eigenvector relation implies, by induction, that
\[
 V^{(z)}Q^n=(-x(z))^nV^{(z)},
 \qquad n\geq 0.
\]
Consequently,
\[
 V^{(z)}P_t
 =
 \sum_{n=0}^{\infty}\frac{t^n}{n!}V^{(z)}Q^n
 =
 \sum_{n=0}^{\infty}\frac{(-x(z)t)^n}{n!}V^{(z)}
 =
 e^{-x(z)t}V^{(z)}.
\]
We shall use this identity for $z=r$ and $z=-r$. In these two cases
$r^\dagger=s$, $(-r)^\dagger=-s$,
and all four geometric sequences are summable because
$0<s<\rho<r<1$.

Since we may write
\[
 m_j=\frac Ar(r^j+\varepsilon(-r)^j),
\]
using $r^j=V_j^{(r)}+s^j$ and $(-r)^j=V_j^{(-r)}+(-s)^j$,
gives the exact decomposition
\[
 m=\frac Ar V^{(r)}
   +\frac{A\varepsilon}{r}V^{(-r)}
   +D,
\]
where the signed remainder $D=(d_j)_{j\geq 1}$ is given by
\[
 d_j=\frac Ar(s^j+\varepsilon(-s)^j).
\]
In particular,
\[
 |d_j|\leq\frac{A(1+\varepsilon)}r s^j.
\]
For the two eigenvector components, we have
\[
 mP_t
 =
 \frac Ar e^{-x(r)t}V^{(r)}
 +\frac{A\varepsilon}{r}e^{-x(-r)t}V^{(-r)}
 +DP_t.
\]
The alternating $V^{(-r)}$ component is exponentially negligible
relative to the positive $V^{(r)}$ component. Indeed,
\[
 x(-r)=a+b+ar+\frac br>\lambda_C,
\]
whereas
\[
 \lambda_C-x(r)
 =ar+\frac br-2\sqrt{ab}>0.
\]
Consequently, $x(-r)>\lambda_C>x(r)$.

Finally, we need to control the remainder. Recall that
$h_j=j\rho^{-j}$ satisfies $P_th=e^{-\lambda_C t}h$.
Since $h_j\geq 1$ for $j\geq 1$ (because $0<\rho<1$),
positivity of the killed semigroup gives
$P_t\mathbf 1\leq P_th=e^{-\lambda_C t}h$.
Therefore
\[
 \|DP_t\|_1
 \leq \sum_{j\geq 1}|d_j|\,\mathbb P_j(T>t)
 \leq e^{-\lambda_C t}\sum_{j\geq 1}|d_j|h_j
 \leq
 \frac{A(1+\varepsilon)}r e^{-\lambda_C t}
 \sum_{j\geq 1}j\left(\frac{s}{\rho}\right)^j.
\]
Because $s<\rho$, the series converges, and hence
\[
 \|DP_t\|_1=O(e^{-\lambda_C t})
           =o(e^{-x(r)t}).
\]
Also,
\[
 \left\|
 \frac{A\varepsilon}{r}
 e^{-x(-r)t}V^{(-r)}
 \right\|_1
 =o(e^{-x(r)t}).
\]
It follows that, in $\ell^1$,
\[
 mP_t
 =
 \frac Ar e^{-x(r)t}V^{(r)}
 +o(e^{-x(r)t}).
\]

The vector $V^{(r)}$ is positive, and
\[
 \|V^{(r)}\|_1
 =\sum_{j\geq 1}(r^j-s^j)
 =\frac{r-s}{(1-r)(1-s)}.
\]
Summing the preceding $\ell^1$ equivalent over $j\geq 1$ gives
\[
 \Pp_m(T>t)
 =
 \frac Ar e^{-x(r)t}
 \frac{r-s}{(1-r)(1-s)}
 +o(e^{-x(r)t}).
\]
Moreover,
\[
 \frac{V_j^{(r)}}{\|V^{(r)}\|_1}
 =
 \frac{(1-r)(1-s)}{r-s}(r^j-s^j)
 =q_j^{(r)}.
\]
Dividing the $\ell^1$ equivalent for $mP_t$ by its total mass
therefore yields
\[
 \left\|
 \mathcal L_m(X_t\mid T>t)-q^{(r)}
 \right\|_{\rm TV}
 \to 0.
\]

Finally,
\[
 \frac{m_{j+1}}{m_j}
 =
 r\frac{1-\varepsilon(-1)^j}{1+\varepsilon(-1)^j}.
\]
Along even and odd values of $j$, this ratio takes the two values
\[
 r\frac{1-\varepsilon}{1+\varepsilon}
 \quad\text{and}\quad
 r\frac{1+\varepsilon}{1-\varepsilon},
\]
respectively. Hence the successive-term ratio does not converge,
even though the conditioned distribution converges to $q^{(r)}$.
\end{proof}

The example shows that an exact domain criterion cannot be formulated
solely through convergence of the raw ratios $m_{n+1}/m_n$. Such a
criterion would have to allow for oscillatory components, such as the
alternating component above, that become asymptotically negligible
under evolution by the killed semigroup because their temporal decay
rates are strictly greater than that of the positive $r$-component.
A smoothed-tail condition, a generating-function condition, or a
condition expressing asymptotic dominance of the positive component
may provide a more intrinsic formulation.

\section{A necessary survival-ratio condition}\label{sec:necessary}

If an initial distribution $m$ lies in the domain of attraction of $q^{(r)}$, then the Markov property gives, for every fixed $u\geq 0$,
\[
 \frac{\mathbb P_m(T>t+u)}{\mathbb P_m(T>t)}
 \to e^{-x(r)u}.
\]
Indeed, by the Markov property, the ratio is the conditional
expectation of $\mathbb P_{X_t}(T>u)$ given survival to time $t$.
This function is bounded by one, so total-variation convergence of the
conditioned laws allows passage to the limiting QSD, under which the
survival probability over an additional time $u$ is $e^{-x(r)u}$.
Hence
\[
 \mathbb P_m(T>t)=e^{-x(r)t}\ell(t),
 \qquad
 \frac{\ell(t+u)}{\ell(t)}\to 1.
\]
This is a genuine necessary condition, not a sufficient one; moreover,
it is a condition on the evolved survival law rather than a complete
description of the initial tail.

\section{Discussion}\label{sec:discussion}

Our results identify substantial subdomains of attraction for every QSD
of the absorbed constant-rate random walk. Theorem~\ref{thm:ratio-tail}
shows that for each $r\in(\rho,1)$ a strictly positive initial
distribution with tail ratio $r$ selects the nonminimal QSD $q^{(r)}$.
This condition covers geometric tails with polynomial, logarithmic or
other ratio-slowly-varying corrections. The oscillatory example of
Section~\ref{sec:counterexample} shows, however, that convergence of the
tail ratio is sufficient but not necessary.

For the minimal QSD, Proposition~\ref{prop:light} shows that
$\sum_{n\geq 1}m_n n\rho^{-n}<\infty$ is sufficient for convergence to
the minimal QSD $q^*$. This includes every initial distribution with
tail ratio strictly less than $\rho$, as well as part of the critical
regime in which the tail ratio is $\rho$. For example, the proposition
applies to $m_n=c\,\rho^n n^{-\beta}$, when $\beta>2$. Thus, for
strictly positive initial distributions having a limiting tail ratio
$r<1$, our results give the following partial classification:
\[
\begin{cases}
 r<\rho &: \text{selection of the minimal QSD }q^*,\\
 r=\rho &: \text{a critical regime not determined by the tail ratio alone},\\
 \rho<r<1 &: \text{selection of the nonminimal QSD }q^{(r)}.
\end{cases}
\]

At the level of generality of the ratio hypothesis, the factor $W(t)$ in
Theorem~\ref{thm:ratio-tail} naturally retains subexponential
information in the sequence $m_n/r^n$; without additional regularity,
one should not expect a universal polynomial prefactor. By contrast,
under the moment condition of Proposition~\ref{prop:light}, the survival
probability has the explicit asymptotic form
\[
 \mathbb P_m(T>t)
 \sim
 \gamma H_m\frac{\rho}{(1-\rho)^2}
 t^{-3/2}e^{-\lambda_C t},
 \qquad
 H_m=\sum_{n\geq 1}m_n n\rho^{-n}.
\]
A complete description of the domains of attraction remains open. There
are two principal gaps. The first is the remainder of the critical
regime: when $m_{n+1}/m_n\to\rho$ but $\sum_{n\geq 1}m_n
n\rho^{-n}=\infty$, the limiting conditional behaviour is not determined
by the present results. The second is the problem of finding an
intrinsic necessary-and-sufficient condition for selecting a given
nonminimal QSD. The counterexample in Section~\ref{sec:counterexample}
shows that convergence of $m_{n+1}/m_n$ to $r$ is not necessary for
selection of $q^{(r)}$. The remaining problem is to find an intrinsic
condition on the initial distribution that is both necessary and
sufficient for this selection. Generating-function singularities and
exponentially smoothed tails are natural candidates for formulating such
a condition.

A natural extension of the results presented here would be to
birth--death processes whose birth and death rates are eventually, or
asymptotically, constant. One expects the limiting rates to determine
the relation between the tail ratio of the initial distribution and the
selected QSD, subject to the spectral effects of the modified rates near
the absorbing boundary. For eventually constant rates, the birth--death
recurrence has geometric solutions outside a finite set, suggesting an
approach based on transfer matrices or finite-rank perturbations of the
constant-rate generator. For asymptotically constant rates, analogous
results would require sufficiently uniform asymptotics for the
recurrence solutions or the associated orthogonal 
polynomials~\cite[Chapter 8]{Anderson1991}. The main
difficulty is that the exact Bessel generating identity used here is
unavailable, so both the exponential scale of the relevant convolution
and the uniform tail bounds needed for total-variation convergence would
have to be obtained by other means.

\bigskip
\noindent
\textbf{Use of AI-assisted technologies.}\quad
The author used Microsoft Copilot, a GPT-based generative AI tool, to
assist with phrasing, exposition, and preliminary checks of intermediate
calculations. The author takes full responsibility for the content of
the paper.

\end{document}